\documentclass[12pt]{amsart}

\usepackage{times,epsf,amssymb,amsmath,hyperref, rlepsf, color}

\begin{document}

\newtheorem{thm}{Theorem}[section]
\newtheorem{lem}[thm]{Lemma}
\newtheorem{cor}[thm]{Corollary}

\theoremstyle{definition}
\newtheorem{defn}[thm]{\bf{Definition}}

\theoremstyle{remark}
\newtheorem{rmk}[thm]{Remark}

\def\square{\hfill${\vcenter{\vbox{\hrule height.4pt \hbox{\vrule width.4pt height7pt \kern7pt \vrule width.4pt} \hrule height.4pt}}}$}

\newenvironment{pf}{{\it Proof:}\quad}{\square \vskip 12pt}

\title{Asymptotic $H$-Plateau Problem in ${\Bbb H}^3$}
\author{Baris Coskunuzer}
\address{MIT, Mathematics Department, Cambridge, MA 02139}
\address{Koc University, Department of Mathematics, Istanbul 34450 Turkey}
\email{bcoskunuzer@ku.edu.tr}
\thanks{The author is supported by Fulbright Grant, and TUBITAK 2219 Grant.}

\maketitle


\newcommand{\Si}{S^2_{\infty}({\Bbb H}^3)}
\newcommand{\PI}{\partial_{\infty}}
\newcommand{\BH}{{\Bbb H}^3}
\newcommand{\BR}{\Bbb R}
\newcommand{\BC}{\Bbb C}
\newcommand{\BZ}{\Bbb Z}

\newcommand{\e}{\epsilon}
\newcommand{\wh}{\widehat}
\newcommand{\p}{\mathcal{P}}
\newcommand{\s}{\Sigma}
\newcommand{\I}{\mathcal{I}}
\newcommand{\F}{\mathcal{F}}
\newcommand{\A}{\mathcal{A}}
\newcommand{\U}{\mathcal{U}}
\newcommand{\B}{\mathbf{B}}

\begin{abstract}

We show that for any Jordan curve $\Gamma$ in $\Si$ with at least one smooth point, there exists an embedded $H$-plane $\p_H$ in $\BH$ with $\PI \p_H =\Gamma$ for any $H\in [0,1)$.
\end{abstract}

\section{Introduction}

There are two versions of the asymptotic Plateau problem. The first version asks the existence of a least area plane $\p$ in $\BH$ asymptotic to a given simple closed $\Gamma$ in $\Si$, i.e. $\PI \p =\Gamma$. In this version, there is a topological restriction on the surface $\p$ to be in the type of disk. The other version asks the existence of an area minimizing surface $\s$ in $\BH$ asymptotic to a given collection of Jordan curves $\wh{\Gamma}$ in $\Si$, i.e. $\PI \s =\wh{\Gamma}$. In the latter version, there is no a priori topological restriction on the surface $\s$, hence $\s$ can have positive genus depending on given $\wh{\Gamma}$. Anderson gave positive answers to both questions 3 decades ago \cite{A1, A2}.

Constant mean curvature (CMC) surfaces are natural generalizations of minimal surfaces, and in many cases, the results related to minimal surfaces are studied to see whether they can be generalized to CMC setting. In our case, we will call this natural generalization as {\em asymptotic $H$-Plateau problem}. A decade after Anderson's result, the second version of the asymptotic Plateau problem was generalized to CMC case by Tonegawa \cite{To}. Tonegawa showed that for any given collection of Jordan curves $\wh{\Gamma}$ in $\Si$, there exists a minimizing $H$-surface $\s_H$ in $\BH$ with $\PI \s_H=\wh{\Gamma}$ where $H\in[0,1)$. Indeed, both Anderson and Tonegawa used geometric measure theory methods, and the solutions are automatically smoothly embedded surfaces by the regularity results of GMT. The survey \cite{C2} gives a fairly complete account of the old and new results on the problem.

On the other hand, for the generalization of the first (plane) version to the CMC case, the only result came out a few years ago by Cuschieri \cite{Cu}. He showed the existence of {\em immersed} $H$-planes asymptotic to given smooth Jordan curve in $\Si$ by using PDE techniques.

In this paper, we give positive answer to the asymptotic $H$-Plateau problem for a larger family of curves. Furthermore, we show that these solutions are indeed \textit{embedded}.

\begin{thm} Let $\Gamma$ be a simple closed curve $\Gamma$ in $\Si$ with at least one $C^1$-smooth point. Then, for any $H\in[0,1)$, there exists a properly embedded $H$-plane $\p_H$ in $\BH$ with $\PI \p_H=\Gamma$.
\end{thm}

Our techniques are also valid for $H=0$ case, and we are able to reprove the existence of least area planes in $\BH$. Hence, with this result, we filled the gap in Anderson's result in \cite{A2}, too (See Remark \ref{gap}).

On the other hand, our proof is indeed for $-1<H<1$ by considering the orientation. Hence, if you forget the sign of the mean curvature, the above theorem shows that for a given Jordan curve $\Gamma$ in $\Si$ and $H\in(0,1)$, there exist {\em a pair of} complete, embedded $H$-planes $\p_H^+$ and $\p_H^-$ with $\PI\p_H^\pm=\Gamma$ (See Remark \ref{pair}).

Recently, Meeks, Tinaglia and the author constructed nonproperly embedded $H$-planes in $\BH$ for any $H\in[0,1)$ \cite{CMT}, where the asymptotic boundary is a pair of infinite lines in $\Si$. Here, we also show that if the Jordan curve $\Gamma$ in $\Si$ is smooth enough, then the minimizing $H$-planes $\p_H$ with $\PI\p_H=\Gamma$ are {\em properly embedded} in $\BH$ (Corollary \ref{proper}).

In the final section, we discuss the questions of the {\em generic uniqueness} of the $H$-planes in $\BH$, and the {\em foliations of} $\BH$ with $H$-planes, and give an outline to solve them.

The organization of the paper is as follows: In the next section, we go over the basic notions, and the related results. In Section 3, we prove the main theorem for least area planes ($H=0$ case). In Section 4, we show the existence of embedded $H$-planes in $\BH$. Finally in Section 5, we give some concluding remarks.

\section{Preliminaries}

In this section, we will go over the related results which will be used in the following sections. For further details of the notions, and the results used in this paper, one can see the survey \cite{C2}.

\begin{defn} Let $M$ be a $3$-manifold.
\begin{itemize}
\item {\em Minimal surface:} A surface $S$ in $M$ is minimal if the mean curvature vanishes everywhere on $S$.

\item {\em Least Area Disk:} A compact disk $D$ in $M$ with $\partial D=\Gamma$ is the least area disk in $M$ if it has the smallest area among the disks in $M$ with the boundary $\Gamma$.

\end{itemize}
\end{defn}

Note that minimal surfaces are the critical points of the area functional. Least area disks and area minimizing surfaces are the minimum of the area functional in the corresponding spaces.

A natural generalization of minimal surfaces are CMC surfaces ($H$-surfaces). They can be defined as the critical points of the area functional with a volume constraint as follows. For an immersion $u:D^2\to M$, the critical points of the following variational problem are immersed disks with constant mean curvature $H$ \cite{Gu}. $$\F_H(u)=\int_{D^2} (|u_x|^2+ |u_y|^2) +\dfrac{4}{3}H [u\cdot(u_x\times u_y)] \ dxdy$$  Here, the second summand in the integral represents the volume constraint.

We can reformulate this variational problem in a different way so that it will be independent of the parametrization of the surface \cite{C5}. Let $\Sigma$ be a surface in $M$ with boundary $\alpha$. We fix a surface $T$ in $M$ with $\partial T = \alpha$, and define $\Omega$ to be the domain bounded by $T$ and $\Sigma$. Again, let $$\I_H(\Sigma)= Area(\Sigma) + 2 H Vol(\Omega)$$ If $\Sigma$ is a critical point of the functional $\I_H$ for any variation $f$, then this will imply $\Sigma$ has constant mean curvature $H$. Note also that critical point of the functional $\I_H$ is independent of the choice of the surface $T$ since if $\widehat{I}_H$ is the functional which is defined with a different surface $\widehat{T}$, then $\I_H - \widehat{\I}_H = C$ for some constant $C$. Note that to keep the solution surface away from $T$, one needs a convexity condition on $T$ (e.g. $H_0$-convex for $H_0>H$) to employ the maximum principle \cite{C5}.

\begin{defn} Let $M$ be a $3$-manifold.
\begin{itemize}
\item {\em $H$-surface:} A surface $S$ in $M$ is $H$-surface if the mean curvature is equal to $H$ everywhere on $S$.


\item {\em Minimizing $H$-disk:} A compact disk $D$ in $M$ with $\partial D=\Gamma$ is a minimizing $H$-disk in $M$ if $\I_H(D)$ (or equivalently $\F_H(D)$) has the smallest value among the disks in $M$ with the boundary $\Gamma$.

\end{itemize}
\end{defn}

\subsection{Embedded solutions to the $H$-Plateau problem:} Here, we will quote the generalization of Meeks-Yau's embeddedness result \cite{MY} to $H$-disks.

\begin{defn} [$H_0$-convex domains] Let $\Omega$ be a compact $3$-manifold with piecewise smooth boundary. We call $\Omega$ an {\em $H_0$-convex domain} if

\begin{itemize}
\item The mean curvature vector $\mathbf{H}$ always points towards inside $\Omega$ along the smooth parts of $\partial \Omega$,

\item  The mean curvature $|\mathbf{H}(p)|\geq H_0$ for any smooth point  $p\in\partial \Omega$,

\item Along the nonsmooth parts of $\partial \Omega$, the inner dihedral angle is less than $\pi$.

\end{itemize}
\end{defn}

The following lemma is on the embeddedness of the solutions of $H$-Plateau problem for $H$-extreme curves.

\begin{lem} \label{embed} \cite{C5} Let $M$ be a compact $H_0$-convex ball. Let $\Gamma$ be a Jordan curve in $\partial M$. Then, for any $H\in[0,H_0)$, there exists a minimizing $H$-disk $\Sigma_H$ in $M$ with $\partial \Sigma_H=\Gamma$, and any such $\Sigma_H$ is embedded.
\end{lem}

The following lemma is know as maximum principle.

\begin{lem} \cite{Gu} \label{maximum} [Maximum Principle] Let $\Sigma_1$ and $\Sigma_2$ be two surfaces in a Riemannian $3$-manifold which intersect at a common point tangentially. If $\Sigma_2$ lies in positive side (mean curvature vector direction) of $\Sigma_1$ around the common point, then $H_1$ is
strictly less than $H_2$ ($H_1 < H_2$) where $H_i$ is the mean curvature of $\Sigma_i$ at the common point.
\end{lem}

\subsection{$H$-Planes in $\BH$:}

Now, we restrict ourselves to $\BH$. Let $\Gamma$ be a simple closed curve in $\Si$. $\Gamma$ separates $\Si$ into two open disks, say $D^+$ and $D^-$. We fix an orientation, and consider the mean curvature $H$ with sign depending on the direction of the mean curvature vector, i.e. $-1<H<1$ such that $\p_H\sim D^- \Rightarrow H\sim +1$ and $\p_H\sim D^+ \Rightarrow H\sim -1$ \cite{C3, To}.

Now, by using the definitions above, we define the least area planes and the minimizing $H$-planes in $\BH$.

\begin{defn} [Least Area Plane] Let $\p$ be a complete surface in $\BH$ which is topologically in the type of a disk. We call $\p$ a least area plane in $\BH$, if any compact subdisk $D$ in $\p$ is a least area disk.
\end{defn}

\begin{defn} [Minimizing $H$-plane] Fix $H\in(-1,1)$. Let $\p_H$ be a complete surface in $\BH$ which is topologically in the type of a disk. We call $\p_H$ a minimizing $H$-plane in $\BH$, if any compact subdisk $D$ in $\p_H$ is a minimizing $H$-disk.
\end{defn}

For a given surface $S$ in $\BH$, we define the {\em asymptotic boundary} of $S$ as follows. If $\overline{\BH}=\BH\cup \Si$ is the natural (geodesic) compactification of $\BH$, and $\overline{S}$ is the closure of $S$ in $\overline{\BH}$, then the asymptotic boundary $\PI S$ of $S$ defined as $\PI S=\overline{S}\cap \Si$.

Now, we define the shifted convex hulls $CH_H(\Gamma)$ as generalizations of the convex hulls in $\BH$ \cite{C2,C3}. Fix $H\in(-1,1)$. Let $\Gamma$ be a Jordan curve in $\Si$. Let $\alpha$ be a  a round circle in $\Si$ with $\alpha\cap \Gamma=\emptyset$. Let $\p^\alpha_H$ be the unique $H$-plane in $\BH$ with $\PI\p^\alpha_H=\alpha$. $\alpha$ separates $\Si$ into two open disks $\Delta_\alpha^+$ and $\Delta_\alpha^-$. Similarly, $\p^\alpha_H$ divides $\BH$ into two domains $\Omega_H^{\alpha +}$ and $\Omega_H^{\alpha -}$ where $\PI \Omega^{\alpha \pm}_H= \Delta_\alpha^\pm$. We will call these regions as \textit{$H$-shifted halfspaces}. If $\Gamma\subset\Delta_\alpha^+$, then we will call $\Omega^{\alpha +}_H$ a \textit{supporting $H$-shifted halfspace}. Similarly, if $\Gamma\subset\Delta_\alpha^-$, then we will call $\Omega^{\alpha -}_H$ a \textit{supporting $H$-shifted halfspace}.

\begin{defn} [Shifted Convex Hull] Let $\Gamma$ be a simple closed curve in $\Si$. Fix $H\in(-1,1)$. Then the {\em $H$-shifted convex hull} of $\Gamma$, $CH_H(\Gamma)$ is defined as the intersection of all supporting closed $H$-shifted halfspaces $\Omega^{\alpha \pm}_H$ of $\BH$. For $H=0$, this is the usual convex hull definition in $\BH$, i.e. $CH(\Gamma)=CH_0(\Gamma)$.
\end{defn}

Now, the generalization of convex hull property of minimal surfaces in $\BH$ to $H$-surfaces in $\BH$ is as follows \cite{C3, To}.

\begin{lem} \label{convex} \cite{To}, \cite{C3}
Let $\Sigma$ be a $H$-surface in $\BH$ where $\PI\Sigma = \Gamma$ and $|H|<1$. Then $\Sigma$ is in the $H$-shifted convex hull of $\Gamma$, i.e. $\Sigma \subset CH_H(\Gamma)$.
\end{lem}

\begin{rmk} \label{convex2} Notice that the result is true for any $H$-surface. This is a straightforward generalization of the convex hull property for minimal surfaces. In particular, if $\Sigma$ is an $H$-surface in $\BH$ with $\PI \Sigma_H=\Gamma$, then $\Sigma$ cannot go into a nonsupporting $H$-shifted halfspaces of $\Gamma$, as we can foliate these halfspaces with $H$-planes, and the first point of touch gives a contradiction with the maximum principle.
\end{rmk}

\section{Existence of Least Area Planes}

In this section, we will focus on $H=0$ case. In other words, we will consider the original asymptotic Plateau problem, and show the existence of smoothly embedded least area planes $\p$ in $\BH$ with $\PI \p =\Gamma$ for a given Jordan curve $\Gamma$ in $\Si$.

Note that Gabai showed the existence of least area planes in $\BH$ in \cite{Ga} by using Hass and Scott's techniques. Recently, Ripoll and Tomi also showed the existence of complete embedded minimal planes in Hadamard manifolds \cite{RT}.

In this paper, we will adapt Anderson's techniques in \cite{A2} to construct minimizing $H$-planes. To generalize his techniques, we need to fill a gap in the proof for least area plane case. Following remark explains the problem.

\begin{rmk} \label{gap} (Gap in \cite[Theorem 4.1]{A2}]) Anderson showed the existence of least area planes in \cite[Theorem 4.1]{A2}. He basically generalized the techniques he used for absolutely area minimizing surface case to the plane case. In particular, let $\Gamma$ be a Jordan curve in $\Si$, and $\{D_n\}$ be a sequence of least area disks in $B_n(0)$ with $\partial D_n=\gamma_n \subset \partial B_n(0)$ where $\gamma_n\to \Gamma$. Then, the idea is to show the existence of a subsequence of $\{D_n\}$ converging to a least area plane $\p$ in $\BH$ with $\PI \p=\Gamma$. In particular, for every fixed compact domain $K$ in $\BH$, he showed the sequence $D_n\cap K=D_n^K$ has a subsequence converging to a smooth disk in $K$ by using the compactness and regularity results of GMT. Then, by using the diagonal sequence argument, he obtained a limit least area plane $\p$ with $\partial \p=\Gamma$.

However, to use the compactness result in this approach, one needs a uniform area bound on the disks $\{D_n^K\}$. Let $\Gamma_n^K$  be the collection of simple closed curves $D_n\cap \partial K$. If $\{\Sigma_n\}$ was a sequence of area minimizing surfaces, then the area of $\partial K$ would be a uniform area bound for $\{\Sigma_n^K\}$ as $\Gamma_n^K$ bounds a surface in $\partial K$ and $\Sigma_n^K$ is absolutely area minimizing surface with boundary $\Gamma_n^K$. Hence, for any $n$, $Area(\Sigma_n^K)\leq Area(\partial K)$.  However, in the disk case $D_n\cap K$ may contain many disks or planar surfaces in $K$, and the area of $\partial K$ cannot give an upper bound for $Area(D_n^K)$. In particular, the estimate (4.2) ($\mathbf{M}(D_i\llcorner\B_r)\leq \frac{1}{2}Area(S(r))$) in \cite[Theorem 4.1]{A2} is not valid in general.

For example, if $D_n^K$ is $2k$ disjoint disks close to the equator disk in $K=\B_R$, and the area of $D_n^K$ would be close to the area of $k$ equator disks, which is much larger than the area of $\partial \B_R$. The main difference with the area minimizing case is that if we  have $k$ annuli $A_1,.. ,A_k$ in $\partial K$ bounding $\Gamma_n^K$, we cannot compare the sum of the areas of the disks with the sum of the areas of the annuli. Because if we replace two disks with an annulus in $D_n$ we would get a genus $1$ surface, which is no longer a disk. So, the area of an annulus cannot be compared with the area of the two disks, because of the restriction of the topology on $\{D_n\}$. Since there is no restriction on the topology of the surface for area minimizing surfaces $Area(\partial K)$ gives a uniform bound, but in the least area disk case, $Area(\partial K)$ does not give a uniform bound for $\{D_n^K\}$.
\end{rmk}

Now, we will show the existence of least area plane in $\BH$ by using Anderson's techniques. In order to get the uniform area bound on $\{D_n^K\}$ for fixed compact set $K$, we will use the ideas in \cite{C1}.

\begin{thm} \label{LAplane} Let $\Gamma$ be a simple closed curve with at least one smooth ($C^1$) point in $\Si$. Then, there exists a properly embedded least area plane $\p$ in $\BH$ with $\PI \p=\Gamma$.
\end{thm}

\begin{pf} We will define a special sequence of least area disks $\{D_n\}$ where their restriction to a compact subset $K$ has a uniform area bound. Then, by following \cite{A2}, we get a least area plane.

{\em Notation and the Setting:} Let $\Gamma$ be a Jordan curve with at least one smooth ($C^1$) point in $\Si$. Let $CH(\Gamma)$ be the convex hull of $\Gamma$ in $\BH$. Near the smooth point $p\in \Gamma$, we can have two sufficiently small and close round circles $\tau^+$ and $\tau^-$ in the opposite sides of $\Gamma $ such that $\tau^+\cup \tau^-$ bounds a least area annulus $\A$ in $\BH$ \cite{Wa}. Hence, $\A$ goes through $CH(\Gamma)$.

Fix a point $O$ in $CH(\Gamma)$. Let $\B_n$ be the closed ball in $\BH$ of radius $n$ with center $O$. Then, for sufficiently large $N_0$, $\A\cap CH(\Gamma)$ is in $\B_{N_0}$. Let $\partial CH(\Gamma)=\partial^+CH(\Gamma)\cup \partial^-CH(\Gamma)$. For $n>N_0$, let $\alpha_n^+=\partial \B_n\cap \partial^+ CH(\Gamma)$ and $\alpha_n^-=\partial \B_n\cap \partial^- CH(\Gamma)$. Then, since $\A$ is a least area annulus, by \cite{MY}, the pair $\alpha_n^+\cup\alpha_n^-$ bounds a least area annulus $\A_n$ in $CH(\Gamma)$.

Let $\Omega_n$ be the compact region which $\A_n$ separates from $CH(\Gamma)$. Let $\gamma_n$ be an essential, smooth, simple closed curve in $\A_n$. Then by construction, $\gamma_n \to \Gamma$ as $\PI CH(\Gamma)=\Gamma$. Let $D_n$ be the least area disk in $\Omega_n$ with $\partial D_n=\gamma_n$ \cite{MY}. By construction, $D_n$ is also a least area disk in $\BH$.

Now, we claim that there exists a subsequence of $\{D_n\}$ which converges to a least area plane $\p$ in $\BH$ with $\PI \p=\Gamma$. We will follow the proof of \cite[Theorem 4.1]{A2}. Hence, if we show the estimate (4.2) \cite[Theorem 4.1]{A2} ($\mathbf{M}(D_n\llcorner\B_r)\leq \frac{1}{2}Area(\partial \B_r)$) is valid for our sequence, we are done.

Now, let $\beta$ be a transversal arc in $CH(\Gamma)$ connecting $\partial^+ CH(\Gamma)$ and $\partial^- CH(\Gamma)$ through the point $O$ (centers of the balls $\B_n$). Let $l$ be the length of $\beta$. Consider the following lemma from \cite{C1}.

Let $D_r$ be a least disk in $\B_r$ with $\partial D_r\subset \partial \B_r\cap CH(\Gamma)$. Then, we call $D_r$ {\em nonseparating with respect to $\Gamma$ in $\B_r$} (say wrt-$\Gamma$), if $\partial D_r$ is not an essential curve in the annulus $\partial \B_r\cap CH(\Gamma)$.

\begin{lem} \label{key} \cite[Lemma 4.1]{C1} Let $\Gamma$ be a Jordan curve in $\Si$ with at least one smooth ($C^1$) point. Let $D_r$ be a least area disk, nonseparating wrt-$\Gamma$ in $\B_r$ for $r>N_0$. Then there is a monotone increasing function $F:[N_0,\infty)\to \BR^+$ such that $F(r)\rightarrow \infty$ as $r \rightarrow \infty$, and $d(O,D_r)>F(r)$ where $d$ is the distance.
\end{lem}

Let $R_0>0$ be sufficiently large so that $F(R_0)>l$. Now, we will prove the uniform bound $|D_n\cap \B_r|<C_r$ for $r>R_0$ where $|.|$ represents the area.  

By \cite[Lemma 4.2]{A2}, for a given least area disk $D_n$, $D_n\cap \B_r$ is a collection of disjoint embedded disks for any generic $r>0$. This is simply because $D_n\cap \B_r$ is a surface for generic $r$, and by the convexity of $\B_r$, any component in $D_n\cap \B_r$ must be a disk.


We claim that for $n>\max\{N_0,R_0\}$, $D_n\cap \B_r$ is just a disk (only one component) for any generic $r>R_0$. Assume on the contrary, and let $E_1, E_2$ be two such disks in $D_n\cap \B_{r_0}$. Since $D_n$ is an embedded disk in $\B_n$, there is a path $\alpha$ connecting $E_1$ and $E_2$ in $D_n$. Let $r_0<r'<n$ be the smallest radius such that $E_1$ and $E_2$ are in the same component $\wh{E}$ of $D_n\cap B_{r'}$. Hence, $\wh{E}$ is a least area disk, nonseparating  wrt-$\Gamma$ in $\B_{r'+\e}$ for sufficiently small $\e>0$. Hence by the Lemma \ref{key}, $d(O, \wh{E})>F(r')$. Hence, we get $d(O, \wh{E})\leq d(O, E_i)\leq l<F(R)<F(r')$ which gives a contradiction.


This proves that for $n>\max\{N_0,R_0\}$ and for any $r>R_0$, $D_n\cap \B_r$ is a disk (only one component). Hence, for any fixed $r>R_0$, $|D_n\cap \B_r|\leq |\partial\B_r|=C_r$, which gives the desired uniform bound. This proves the estimate (4.2) of \cite[Theorem 4.1]{A2} is valid for our sequence $\{D_n\}$. The proof follows.
\end{pf}

\section{Minimizing $H$-Planes}

In this section, we will show the existence of the solutions of asymptotic $H$-Plateau problem in $\BH$ for $H\in(-1,1)$. In particular, we will generalize the techniques in the previous section to the CMC case, and show that for any simple closed curve $\Gamma$ in $\Si$ with one smooth point, there exists an embedded $H$-plane $\p_H$ with $\PI \p_H=\Gamma$. First, we need to generalize Lemma \ref{key} proven in \cite{C1} to the minimizing $H$-planes in $\BH$.

We will use the same notation, i.e. let $O, \B_r,\A_r, N_0$ be as in the previous section. Fix $H\in(-1,1)$. Let  $\wh{A}_\Gamma^H=\partial \B_r\cap CH_H(\Gamma)$ be an annulus in $\partial \B_r$. Again, let $D_r$ be a minimizing $H$-disk in $\B_r$ with $\partial D_r\subset \wh{A}_\Gamma^H$ . Then, we call $D_r$ {\em nonseparating wrt-$\Gamma$ in $\B_r$}, if $\partial D_r$ is \underline{not} an essential curve in the annulus $\wh{A}_\Gamma^H$.

\begin{lem} \label{key2} Fix $H\in(-1,1)$. Let $\Gamma$ be a Jordan curve with at least one smooth ($C^1$) point in $\Si$. For $r>N_0$, let $D_r$ be a minimizing $H$-disk, and nonseparating wrt-$\Gamma$ in $\B_r$. Then there is a monotone increasing function $F:[N_0,\infty)\to \BR^+$ such that $F(r)\rightarrow \infty$ as $r \rightarrow \infty$, and $d(O,D_r)>F(r)$ where $d$ is the distance.
\end{lem}

\begin{figure}[t]

\relabelbox  {\epsfxsize=3.7in

\centerline{\epsfbox{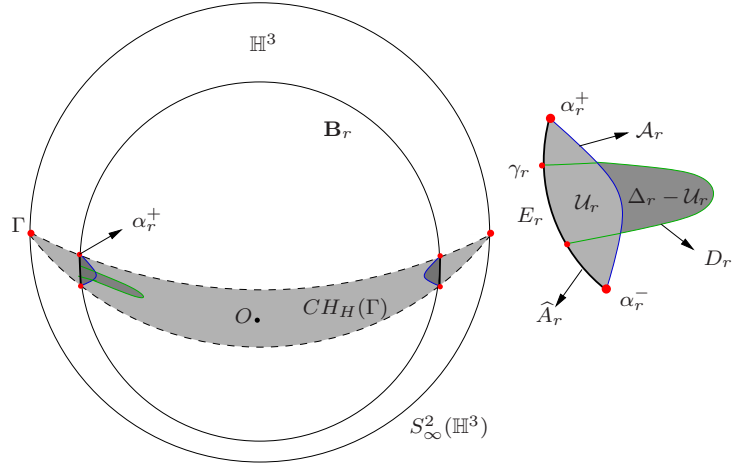}}}

\relabel{1}{\tiny $\Gamma$}
\relabel{2}{\tiny $O$}
\relabel{3}{\tiny $CH_H(\Gamma)$}
\relabel{4}{\tiny $\B_r$}
\relabel{5}{\tiny $\Si$}
\relabel{6}{\tiny $\alpha^+_r$}
\relabel{7}{\tiny $\wh{A}_r$}
\relabel{8}{\tiny $D_r$}
\relabel{9}{\tiny $\gamma_r$}
\relabel{10}{\tiny $\U_r$}
\relabel{11}{\tiny $\Delta_r-\U_r$}
\relabel{12}{\tiny $\alpha^-_r$}
\relabel{13}{\tiny $E_r$}
\relabel{14}{\tiny $\A_r$}
\relabel{15}{\tiny $\alpha^+_r$}
\relabel{16}{\tiny $\BH$}
\endrelabelbox

\caption{\label{key3} \footnotesize If $\gamma_r$ is a nonessential curve in $\wh{A}_r$, the nonseparating minimizing $H$-disk $D_r$ in $\B_r$ with $\partial D_r=\gamma_r$ must belong to $\U_r$. This shows that for any such $D_r$ in $\B_r$, $d(O,D_r)>d(O,\A_r)=F(r)$. }

\end{figure}

\begin{pf} We will adapt the proof of Lemma \ref{key} to this case. Lemma \ref{key} finishes the $H=0$ case. Hence, we can take $H>0$. For $H<0$, the same proof works by changing the orientation. Fix $H>0$.

Recall that $\A_r$ is the least area annulus in $\B_r$ with $\partial \A_r=\alpha_r^+\cup \alpha^-_r$ for generic $r>N$. Let $\wh{A}_r$ be the annulus in $\partial \B_r$ with $\partial \wh{A}_r=\alpha_r^+\cup\alpha^-$. Let $F(r)=d(O,\A_r)$.

We claim that for any minimizing $H$-disk $D_r$ which is nonseparating wrt-$\Gamma$ in $\B_r$, $D_r\cap \A_r=\emptyset$, i.e. $D_r$ stays in the solid torus $\U_r$ in $\B_r$ with $\partial \U_r=\A_r\cup \wh{A}_r$; See Figure \ref{key3}. In particular, this shows that $d(O,D_r)>F(r)$, and the proof will follow.

Now, let $D_r$ be a minimizing $H$-disk and nonseparating wrt-$\Gamma$ in $\B_r$ with $\partial D_r= \gamma_r\subset \wh{A}_r$. Since $D_r$ is nonseparating wrt-$\Gamma$, $\gamma_r$ is a not an essential curve in $\wh{A}_r$. In other words, $\gamma_r$ bounds a disk $E_r$ in $\wh{A}_r$; see Figure \ref{key3}.

Let $\Delta_r$ be the region in $\B_r$ with $\partial \Delta_r=D_r\cup E_r$. Since $D_r$ is a minimizing $H$-disk in $\B_r$ with $H>0$, $I_H(D_r)=|D_r|+2H|\Delta_r|$ is the smallest among all the disks in $\B_r$ with boundary $\gamma_r$. Here, $|.|$ represents the area or the volume of the corresponding region.

Assume that $d(O,D_r)<d(O,\A_r)$. Recall that $\U_r$ is the region in $\B_r$ with $\partial \U_r=\A_r\cup \wh{A}_r$. Let $\Delta'_r=\Delta_r\cap \U_r$. Then as $d(O,D_r)<d(O,\A_r)$, $|\Delta'_r|< |\Delta_r|$. Furthermore, let $T'=\A_r\cap\partial \Delta'_r$ be the planar region in $\A_r$ with $\partial T'=\beta$. Let $T$ be the planar region in $D_r$ with $\partial T= \beta$. Since $D_r$ is a nonseparating disk, and $\A_r$ is an annulus, $D'=(D_r-T)\cup T'$ is also a disk in $\B_r$ with $\partial D'=\gamma_r$. Furthermore, as $\A_r$ is a least area annulus, $|T'|<|T|$ and hence $|D'|<|D_r|$. This implies $\I_H(D')=|D'|+2H|\Delta'_r|<|D_r|+2H|\Delta'_r|=\I_H(D_r)$. However, $D_r$ is a minimizing $H$-disk in $\B_r$. This is a contradiction. This proves $D_r\cap \A_r=\emptyset$, and $D_r\subset \U_r$. Hence, this shows that $d(O,D_r)>d(O,\A_r)=F(r)$. The lemma follows.
\end{pf}

Now, we prove the main result of the paper.

\begin{thm} \label{H-plane} [Existence of $H$-Planes] Let $\Gamma$ be a Jordan curve in $\Si$ with at least one smooth ($C^1$) point. Let $H\in(-1,1)$. Then, there exists a properly embedded minimizing $H$-plane $\p_H$ in $\BH$ with $\PI \p_H=\Gamma$.
\end{thm}

\begin{pf} Fix $H\in(-1,1)$. Given $\Gamma$ in $\Si$, we will use the same setup as before, i.e. let $\A,\A_r,\wh{A}_r,O,\B_r, N_0$ be as in the previous lemma.



Let $\Omega_n=\B_n\cap CH_H(\Gamma)$ for $n>N_0$. Let $\gamma_n$ be an essential smooth curve in $\wh{A}_n\subset\partial \Omega_n$. Then by construction, $\gamma_n \to \Gamma$ as $\PI CH_H(\Gamma)=\Gamma$. Recall that the mean curvature of the geodesic sphere of radius $R$ is $\coth{R}$. Hence, $\B_n$ is $1$-convex for any $n$. Then, by Lemma \ref{embed}, there exists an embedded minimizing $H$-disk $D_n$ in $\B_n$ with $\partial D_n=\gamma_n$. Here, the sign of $H$ determines the minimizing $H$-disk as there are two minimizing $H$-disks in $\B_n$ facing each other for $H\neq 0$ \cite{C5}. Furthermore, since $\gamma_n\subset CH_H(\Gamma)$, by Lemma \ref{convex} (see also Remark \ref{convex2}), $D_n\subset CH_H(\Gamma)$. Hence, this implies $D_n\subset \Omega_n=\B_n\cap CH_H(\Gamma)$.

We claim that the sequence of embedded minimizing $H$-disks $\{D_n\}$ converges to a minimizing $H$-plane $\p_H$. First, we claim that the sequence of minimizing $H$-disks $\{D_n\}$ has a convergent subsequence in any compact set $K$ in $\BH$. Then, by using the diagonal sequence argument as before, we will get a limit minimizing $H$-plane in $\BH$.

Consider the closed balls $\B_{K}$ with center $O$. By using convexity of $\B_K$, for generic $K>N_0$, we can assume $D_n\cap \B_K$ is a collection of disks for any $n$ by \cite[Lemma 4.2]{A2}.

Now, let $\beta$ be a transversal arc in $CH_H(\Gamma)$ connecting $\partial^+ CH_H(\Gamma)$ and $\partial^- CH_H(\Gamma)$ through the point $O$. Let $l$ be the length of $\beta$. Fix a generic $K_0>N_0$ such that $F(K_0)>l$. Then by Lemma \ref{key2}, and the proof of Theorem \ref{LAplane}, $D^0_n=D_n\cap \B_{K_0}$ is a closed disk (only one component) for any $n$. Hence, the sequence of minimizing $H$-disks $\{D^0_n\}$ has a uniform area bound, say $|D^0_n|<|\partial \B_{K_0}|$. Then, by following the proof of \cite[Theorem 4.1]{A2}, the existence of smoothly embedded minimizing $H$-plane $\Sigma_H$ can be shown as follows:

With the uniform area bound, the sequence of minimizing $H$-disks $\{D_n^0\}$ has a convergent subsequence in $B_{K_0}$ by the compactness theorem for integral currents.  Hence by considering them as integral currents, a subsequence of $\{D_n^0\}$ converges to a properly embedded minimizing $H$-disk $D^0$ in $\B_{K_0}$. Let $K_i$ be a monotone increasing sequence with $K_i\nearrow \infty $. For $K_1>K_0$, by starting with this subsequence, get another subsequence converging on $\B_{K_1}$. By iterating this process and diagonal sequence argument, we get a sequence of integral currents $\{D_n\}$ converges on compacts to the integral current $\Sigma_H$ in $\BH$. By Allard's regularity, the convergence is smooth on compact sets. Also, the asymptotic boundary of the support of $\Sigma_H$ is $\Gamma$ by the convex hull property (Lemma \ref{convex}), i.e. $\PI \Sigma_H=\Gamma$ as $\PI CH_H(\Gamma)=\Gamma$.

The limit of minimizing $H$-disks $\Sigma_H$ is a minimizing $H$-surface. Hence, by \cite{Gu}, for any point $p$ in the support of $\Sigma_H$, there exists $\e>0$ with $\B_\e(p)\cap \Sigma_H$ is a smooth embedded disk. Hence, the support of $\Sigma_H$ is smoothly embedded surface. Finally, since the convergence is smooth in compact sets, and $\{D_n\}$ is a sequence of embedded disks, $\Sigma_H$ is a complete minimizing $H$-plane in $\BH$ with $\PI\Sigma_H=\Gamma$. The proof follows.
\end{pf}

\begin{rmk} \label{pair} [Pair of $H$-planes] Notice that if you forget the sign of the mean curvature $H$, the theorem above shows that for a given Jordan curve $\Gamma$ in $\Si$, there exist two minimizing $H$-planes $\p_H^+$ and $\p_H^-$ with $\p_H^\pm=\Gamma$ for $H\in(0,1)$. Furthermore, these $H$-planes are disjoint, and the convex sides are facing each other.
\end{rmk}

\begin{rmk} The above result was shown for special families of curves in $\Si$ like star-shaped curves \cite{GS}, "mean convex" curves \cite{NS}, and graph over a line \cite{RT2}, where they showed that the area minimizing surface $\Sigma_H$ is indeed a graph over a geodesic subspace in $\BH$.
\end{rmk}

\begin{cor} \label{proper} [Properly Embeddedness] Let $\Gamma$ be a $C^{3,\alpha}$ smooth Jordan curve in $\Si$. Let $\p_H$ be a minimizing $H$-plane in $\BH$ with $\PI \p_H=\Gamma$. Then, $\p_H$ is properly embedded in $\BH$.
\end{cor}

\begin{pf} By \cite{To}, since $\Gamma$ is $C^{3,\alpha}$ smooth, $\p_H$ is regular near infinity. In particular, there exists a $\rho>0$ such that in the upper half space model of $\BH$, $\p_H\cap \{z<\rho\}$ is a graph over $\Gamma\times (0,\rho)$. Then, for sufficiently large $N$, $\partial\B_N\cap \p_H$ is a Jordan curve $\gamma$ in $\partial\B_N$. Hence, $\B_N\cap \p_H=\mathcal{D}$ is a minimizing $H$-disk in $\B_N$ with $\partial \mathcal{D}=\gamma$ by the definition of minimizing $H$-plane. By Lemma \ref{embed}, $\mathcal{D}$ is properly embedded. Since $\p_H-\mathcal{D}$ is a graph over $\Gamma\times(0,\rho)$ in the upper half space model \cite{To}, the proof follows.
\end{pf}

\begin{rmk} Indeed, it can be showed that the above result is true for far more generality. By using the techniques in \cite{C1}, one can naturally generalize the above result to Jordan curves in $\Si$ with at least one smooth ($C^1$) point. Even though there is no regularity near infinity in that case, by using Lemma \ref{key2}, the arguments in \cite{C1} can easily be adapted.

Note also that, Meeks, Tinaglia and the author recently showed that there exists a {\em nonproperly} embedded complete $H$-plane $\Sigma_H$ in $\BH$  for any $H\in[0,1)$ \cite{CMT}. In particular, $\Sigma_H$ is an $H$-plane between two rotationally invariant $H$-catenoids $\mathcal{C}_1$ and $\mathcal{C}_2$ where $\Sigma_H$ spirals into $\mathcal{C}_1$ in one end, and spirals into $\mathcal{C}_2$ in the other end. Hence, $\PI \Sigma_H$ is a pair of infinite lines $l^+$ and $l^-$ in $\Si$. Here, if $\PI \mathcal{C}_i=\alpha_i^+\cup\alpha_i^-$, and $A^\pm$ is the annuli in $\Si$ with $\partial A^\pm_i=\alpha_1^\pm\cup\alpha_2^\pm$, then $l^+\subset A^+$ and $l^-\subset A^-$ where $l^\pm$ spirals into $\alpha_1^\pm$ in one end, and spirals into $\alpha_2^\pm$ in the other end. 
\end{rmk}

\section{Final Remarks}

\subsection{Generic Uniqueness of minimizing $H$-Planes} \

The generic uniqueness results for minimizing $H$-surfaces in $\BH$ \cite{C3} can naturally be generalized to our context, i.e minimizing $H$-planes. In particular, for fixed $H\in(-1,1)$, let $\Sigma_1$ and $\Sigma_2$ be minimizing $H$-planes with $\PI \Sigma_i=\Gamma_i$ where $\Gamma_1$ and $\Gamma_2$ are disjoint simple closed curves in $\Si$. Then, by using Meeks-Yau exchange roundoff trick, it can be showed $\Sigma_1$ and $\Sigma_2$ are disjoint, too. By using this, and similar ideas to \cite{C3}, it can be showed that any simple closed curve $\Gamma$ bounds either a unique minimizing $H$-plane $\Sigma$ in $\BH$, or there are two canonical disjoint minimizing $H$-planes $\Sigma^+$ and $\Sigma^-$ with $\PI \Sigma^\pm=\Gamma$. Hence, foliating an annular neighborhood of $\Gamma$ in $\Si$ with simple closed curves, and considering the canonical $H$-planes constructed, one can get generic uniqueness result as in \cite{C3}. In particular, this shows that {\em for fixed $H\in(-1,1)$, a generic Jordan curve $\Gamma$ in $\Si$ bounds a unique minimizing $H$-plane $\Sigma$ in $\BH$ with $\PI\Sigma=\Gamma$}.

\subsection{Foliations of $\BH$ by $H$-planes} \

Similar to the previous part, it is also possible to show that two $H$-planes in $\BH$ with different $H$ values, which are asymptotic to the same asymptotic curve in $\Si$ are disjoint. In particular, for a given Jordan curve $\Gamma$ in $\Si$, and given $-1<H_1<H_2<1$, it can be showed that the minimizing $H$-planes $\p_{H_1}$ and $\p_{H_2}$ in $\BH$ with $\p_{H_i}=\Gamma$ are disjoint by using the ideas in \cite{C4}. Hence, if $\Gamma$ bounds a unique minimizing $H$-plane for any $H\in(-1,1)$, it can be showed that the family of planes $\mathcal{F}_\Gamma=\{\p_H \ | \ \PI\p_H=\Gamma \mbox{ and } -1<H<1\}$ foliates $\BH$. By the disjointness of the planes $\{\p_H\}$ for different $H$'s, in order to get the foliation, all one needs to show is that there is no gap between the planes $\{\p_H\}$ by using the uniqueness. By the arguments in \cite{C4}, a gap between the planes, say between $\{\p_H \ | \ -1<H\leq H_0\}$ and $\{\p_H \ | \ H_0<H<1\}$, implies the nonuniqueness for $H_0$-planes with $\p_{H_0}=\Gamma$. This gives a contradiction, and shows that $\mathcal{F}_\Gamma$ foliates $\BH$. For example, if $\Gamma$ is a star-shaped curve in $\Si$, then the family $\mathcal{F}_\Gamma$ foliates $\BH$ \cite{C4, GS, RT2}.


\end{document}